\documentclass[11pt]{amsart}
\usepackage{amsxtra}
\usepackage{amssymb}
\usepackage{graphicx}
\usepackage{showlabels}
\theoremstyle{plain}
\newmuskip\pFqmuskip

\newcommand*\pFq[6][8]{%
	\begingroup % only local assignments
	\pFqmuskip=#1mu\relax
	\mathcode`\,=\string"8000
	\begingroup\lccode`\~=`\,
	\lowercase{\endgroup\let~}\pFqcomma
	{}_{#2}F_{#3}{\left[\genfrac..{0pt}{}{#4}{#5};#6\right]}%
	\endgroup
}
\newcommand{\pFqcomma}{\mskip\pFqmuskip}

\newtheorem{theorem}{Theorem}[section]
\newtheorem{lemma}[theorem]{Lemma}
\newtheorem{proposition}[theorem]{Proposition}
\newtheorem{corollary}[theorem]{Corollary}
\newtheorem{definition}[theorem]{Definition}
\newtheorem{example}[theorem]{Example}

\newtheorem{remark}[theorem]{Remark}
\newtheorem{conjecture}[theorem]{Conjecture}

\newtheorem{question}[theorem]{Question}
\newcommand \bth[1] { \begin{theorem}\label{t#1} }
	\newcommand \ble[1] { \begin{lemma}\label{l#1} }
		\newcommand \bpr[1] { \begin{proposition}\label{p#1} }
			\newcommand \bco[1] { \begin{corollary}\label{c#1} }
		\newcommand \bde[1] { \begin{definition}\label{d#1}\rm }
					\newcommand \bex[1] { \begin{example}\label{e#1}\rm }
						\newcommand \bre[1] { \begin{remark}\label{r#1}\rm }
							\newcommand \bcon[1] { \begin{conjecture}\label{con#1}\rm }
								\newcommand \bque[1] { \begin{question}\label{que#1}\rm }
			\newcommand{\beq}{\begin{equation}}
			\newcommand{\eeq}{\end{equation}}
			\newcommand{\beqa}{\begin{eqnarray}}
			\newcommand{\eeqa}{\end{eqnarray}}
			\newcommand{\beaa}{\begin{eqnarray*}}
				\newcommand{\eaa}{\end{eqnarray*}}
		
								\newcommand {\ele} { \end{lemma} }
									\newcommand {\ethe} {\end{theorem} }
		\newcommand {\epr} { \end{proposition} }
						\newcommand {\eco} { \end{corollary} }
					\newcommand {\ede} { \end{definition} }
				\newcommand {\eex} { \end{example} }
			\newcommand {\ere} { \end{remark} }
		\newcommand {\econ} { \end{conjecture} }
	\newcommand {\eque} { \end{question} }
\newcommand \thref[1]{Theorem \ref{t#1}}
\newcommand \leref[1]{Lemma \ref{l#1}}

\newcommand \coref[1]{Corollary \ref{c#1}}

\def \B {{\mathcal B}}

\def \M {{\mathcal M}}

\def \Cset {{\mathbb C}}
\def \Zset {{\mathbb Z}}
\def \Nset {{\mathbb N}}

\def \ad { {\mathrm{ad}} }

\def \Cset {{\mathbb C}}

\def \Zset {{\mathbb Z}}

\def \Nset {{\mathbb N}}

\def \B  {{\mathcal{B}}}               %mathcal

\def \la {\lambda}
\def \La {\Lambda}

\def \Om {\Omega}

\def \mt  {\mapsto}
\newmuskip\pFqmuskip

\usepackage{hyperref}

\begin{document}
%%%%%%%%%%%%%%%%%%%%%%%%%%%%%%%%%%%%%%%%%%%%%%%%%%%%%%%%%%%%%%%%%%%%%%%%%%%
%%%%%%%%%%%%%%%%%%%%%%    Title    %%%%%%%%%%%%%%%%%%%%%%%%%%%%%%%%%%%%%%%%
\title[Automorphisms of algebras and discrete VOP]
 {Automorphisms of algebras and Bochner`s property for discrete  vector orthogonal  polynomials}
 	
 	\author[E.~Horozov]{Emil Horozov}
 	\address{
 		Institute of Mathematics and Informatics, \\ 
 	 		Bulg. Acad. of Sci., Acad. G. Bonchev Str., Block 8, 1113 Sofia,
 		Bulgaria  	}
 	\email{horozov@fmi.uni-sofia.bg}

\date{\today}

\date{}

\begin{abstract} 
We construct new families of discrete vector orthogonal polynomials that have the property to be    eigenfunctions of some difference operator. They are extensions of Charlier, Meixner and Kravchuk polynomial systems.  The ideas behind our approach lie in the studies of bispectral operators. We exploit automorphisms of associative algebras which transform elementary (vector) orthogonal polynomial systems which are eigenfunctions of a difference operator into other systems of this type. While the extension of Charlier polynomilas is well known it is obtained by different methods. The extension of Meixner polynomial system is new.
\end{abstract}

\maketitle
%%%%%%%%%%%%%%%%%%   Introduction   %%%%%%%%%%%%%%%%%%%%%%%%%%%%%%%%%

\medskip

\section{Introduction}\label{intro}

S. Bochner   \cite{Bo} has classified all systems of orthogonal (with respect to some measure on the real line) polynomials $P_n(x), \; \; n=0, \ldots,$    that are also  eigenfunctions of a second order differential operator

\begin{equation}
L(x, \partial_x) = A(x)\partial_x^2 + B(x)\partial_x + C(x)  \label{BP1}
\end{equation}
with eigenvalues $\lambda_n$.
Here the coefficients $A, B, C$ of the differential equation do not depend on the index $n$.

A similar problem was solved by O. Lancaster \cite{Lan} and P. Lesky \cite{Les}, although earlier E. Hildebrandt \cite{Hil} has found all needed components of the proof. For more information see the excellent review article by W. Al-Salam \cite{AlS}.

The orthogonality condition, due to a classical theorem by Favard-Shohat is equivalent to  the well known 3-terms recursion relation

\begin{equation}
xP_n = P_{n+1} + \beta(n)P_n + \gamma(n) P_{n-1},  \label{BP2}
\end{equation}
where $\beta(n), \gamma(n)$ are constants, depending on $n$. Here we use  the polynomials normalized by the condition that their highest order coefficient is $1$. The statement of Lancaster's theorem is that all polynomial systems   with such properties are the discrete orthogonal polynomials of Hahn, Meixner, Charlier and Kravchuk.

Here when speaking about  orthogonality we mean orthogonality with respect to a nondegenerate functional, which does not need to be positive definite.

In recent times there is much activity in generalizations and versions of this classical result of Bochner. 

The first one was the generalization by H. L. Krall \cite{Kr}. He classified all  order 4 differential  operators which have a family of orthogonal polynomials as eigenfunctions.
Later many authors found new families of orthogonal polynomials that are eigenfunctions of a differential operator. 

The classical discrete orthogonal polynomials have also  been studied with respect to extending the theorem of Hildebrandt-Lancaster-Lesky. In particular A. Dur\'an and M. de la Iglesia \cite{DdlI} have obtained extensions of  the classical polynomial systems of Hahn, Meixner and Charlier.

An important  role in some of these generalizations plays the ideology of the bispectral problem which was initiated in \cite{DG}. Even translating the Bochner and Krall results into this language already gives a good basis to continue investigations. We formulate it for the case of discrete orthogonal polymials. Let us introduce the function  $\psi(x, n) = P_n(x)$. By $D$ we denote the shift operator in $x$:  $Df(x)= f(x+1)$.    If we write the right-hand side of the 3-term recursion relation as a difference operator $\Lambda(n)$ acting on the discrete variable $n$  then the  3-term recursion relation can be written as

\[
\Lambda(n) \psi(x, n)  = x\psi(x, n) .
\]
This means that  is an eigenfunction of the discrete operator $\Lambda$ with  eigenvalue $x$ and of a difference operator  $L(x, \partial_x)$ with eigenvalue 
$\lambda_n$. In the discrete case instead of the  differential operator we have to 
take a difference operator $L(x, D)$, where     Hence we can formulate the discrete version of Bochner-Krall problem as

{\it Find all systems of orthogonal  polynomials $P_n(x)$   (with respect to some functional $u$) which are eigenvalues of a   difference operator.}

% % % % % % % % % % % % % % % % % % % % % % % % % % % % % %

We also use some  ideas found in the studies of bispectral operators. Before explaining the main results let us introduce one more concept which is central for the present paper. This is the notion of vector orthogonal polynomials (VOP), introduced by J. van Iseghem \cite{VIs}.   Let $\{P_n(x)\}$ be a family of monic polynomials such that $\deg P_n = n$. Assume that they satisfy a $d+2$-term recursion relation, $d \geq 1$ 

\beq
xP_n(x) = P_{n+1} + \sum_{j=0}^{d}\gamma_j (n)P_{n-j}(x) \label{d-ort}
\eeq
with constants  (independent of $x$)  $\gamma_j(n)$, $\gamma_j(d) \neq 0$. Then by a theorem of P. Maroni \cite{Ma} there exist $d$  functionals $u_j, j =0,\ldots,d-1$    on the space of all polynomials  $\Cset[x]$   such that

\[
\begin{cases}
u_k (P_mP_n)    = 0, m > nd + k, n \geq 0,  \\
u_k(P_nP_{n(d+1)+k})    \neq 0, n   \geq 0,    
\end{cases}
\]
for each $k \in N_{d+1} := \{0, \ldots, d-1 \}$. When $d = 1$ this is the  notion of orthogonal polynomials. 

In the last 20-30 years there is much activity in the study of vector orthogonal polynomials and the broader class of multiple orthogonal polynomials.

% % % % % % % % % % % % % % % %

Applications of the d-orthogonal polynomials include the simultaneous Pad\'e
approximation problem where the   multiple orthogonal polynomials appear  \cite{ApKu}.    Also   multiple orthogonal polynomials play important role in random matrix theory   \cite{ApKu, BDK}.   The d-orthogonal polynomials can be obtained from general multiple orthogonal polynomials under some restrictions upon their parameters.

One problem that deserves attention is to find analogs of classical orthogonal polynomials. Several authors  \cite{ACVA, vAssc} have found multiple    orthogonal polynomials, that share a number of properties  with the classical orthogonal polynomials - they have a raising operators, Rodrigues type formulas, Pearson equations for the weights, etc. However one of the features of the classical orthogonal polynomials - a differential operator for which the polynomials are eigenfunction is missing. Sometimes this property is relaxed to the property that the polynomials satisfy linear differential equation, whose coefficients may depend on the index of the polynomial.

  We are looking for polynomials $P_n(x), \; n=0, 1, \ldots$ that are eigenfunctions of a difference  operator $L(x, D)$ with eigenvalues  depending on the  variable $n$ (the index) and which at the same time are eigenfunctions of a difference operator in $n$, i.e. finite-term recursion relation with an eigenfunction, depending only on the  variable $x$.

Our main result includes  an extension  of Meixner polynomials. We construct systems of vector orthogonal polynomials $\{P_n(x)\}$ which are eigenfunctions of a difference operator. It is different from the  family found in \cite{GVZ} and \cite{BDK} except for the first member.
Our approach uses ideas of the bispectral theory from \cite{BHY} but does not use Darboux transformations, which is usually the case, see e.g. \cite{GHH, GY, DdlI}. We use methods introduced in \cite{BHY} which we describe below. Also a well known extension of Charlier polynomials (see \cite{BCZ, VZh}) is presented. The reason to repeat it  is that our construction is a new one in comparison to the techniques of \cite{BCZ, VZh}. However there are some similarities with \cite{VZh}. The authors also use automorphisms of algebras and make a beautiful connection  with representation theory. Our construction is simpler and quite straightforward. The same method was recently  applied to extensions of Hermite and Laguerre polynomials as well as to a family that has no classical analog \cite{Ho}. The present paper could be considered a continuation of \cite{Ho} but  could be read independently.

The method from the present paper and \cite{Ho} can be applied to various versions of vector orthogonal polynomials - matrix, multivariate, etc. This will be done elsewhere.

\medskip
\noindent 

{\it Acknowledgements.} 
The author is deeply grateful to Boris Shapiro for showing and discussing some examples of systems of polynomials studied here and in particular the examples from \cite{ST}. They helped me to guess that the methods from \cite{BHY}  can be useful for the study of vector orthogonal polynomials.    The author is    grateful to the Mathematics Department   of Stockholm University for the hospitality in  April 2015.

Last but not least I am extremely grateful to Prof. T. Tanev, Prof. K. Kostadinov,  and Mrs. Z Karova from the Bulgarian Ministry of Education and Science  and   Prof. P. Dolashka, BAS who  helped me in the difficult situation when I was  sacked  by Sofia university in violations of the Bulgarian laws. This was a retaliation for my attempt to reveal a large -scale  corruption, that involves   highest university and science officials in Bulgaria.\footnote{See, e.g. EMS NEWSLETTER,  http://www.ems-ph.org/journals/newsletter/pdf/2015-12-98.pdf}

\section{Elements of bispectral theory}

%%%%%%%%%%%%%%%%%%%%%%%%%%%%%%%%%%%%%%%%%%%%%%%%%%%%%%%%%%%%%%%%%%%%%%%%%%%%%%

The following introductory material is mainly borrowed from \cite{BHY}. 
Below we present the  difference-difference  version of the general bispectral problem which is  suitable in the set-up of discrete orthogonal polynomial sequences. 

%%%%%%%%%%%%%%%%%%%%%%%%%%%%%%%%%%%%%%%%%%%%%%%%

For $i=1,2$, let $\Om_i$ be two open subsets of $\Cset$   such that $\Om_1$ is invariant 
under the translation operator
\[
D \colon x \mt x+1, \; x \in \Om_1
\]
and its inverse  $D^{-1}$, while  $\Om_2$ is invariant 
under the translation operator
\[
T \colon n \mt n+1
\]
and its inverse  $T^{-1}$.

A complex analytic difference operator on $\Om_1$ 
is a finite sum of the form
\[
\sum_{k \in \Zset} c_k(x) D^k, 
\]
where $c_k \colon \Om_1 \to \Cset$ are analytic functions.   In the same way we define   complex analytic difference operators on $\Om_2$ to be 
 finite sums of the form
\[
\sum_{k \in \Zset} s_k(n) T^k, 
\]
where $s_k \colon \Om_1 \to \Cset$ are analytic functions.

By $\B_1$ we denote an  algebra with unit, consisting of difference  operators     $L(x, D)$  in one variable $x$.  By 
 $\B_2$ we denote an algebra of difference operators $\Lambda(n, T)$.  Denote by $\M$ the space of complex analytic functions 
on $\Om_1 \times \Om_2$.    The space $\M$ is naturally equipped with the structure of bimodule
over the algebra  of  analytic difference operators $L(x, D)$ on $\Om_1$ and the difference operators   $\La(n, T)$   on $\Om_2$.

Assume that there exists an algebra isomorphism $b \colon \B_1 \to \B_2$ and 
an element $\psi \in \M$ such that 
\[
P \psi = b(P)\psi , \quad \forall P \in \B_1.
\]

A  discrete-discrete bispectral function is by definition an element  $\M$ (i.e., 
an  analytic function) 
\[
\psi \colon \Om_1 \times \Om_2 \to \Cset
\]
for which there exist   analytic difference    operator  
$L(x, D)$ and $\La(n, T)$ on 
$\Om_1$ and $\Om_2$, and   analytic functions 
\[
\theta(x)   \quad \text{and } \quad
\la(n), 
\]
such that

\begin{equation}
\begin{split}
\label{bisp3}
L(x, D) \psi(x,n) &=  f(n) \psi(x,n) \\
 \La(n, T_n) \psi(x, n) &=  \theta(x) \psi(x,n) 
\end{split}
\end{equation}
on $\Om_1 \times \Om_2$. Let us introduce the subalgebras $K_i \,\; i= 1,2$ of $\B_i$ to be the algebras of functions in $x$ (respectively in $n$). In fact, as we would be interested in VOP, we will consider only the case  when $\theta(x) \equiv x$.  We will assume that $\psi(x,n)$ is a nonsplit function of $x$ and $n$ in the sense that it satisfies the condition

$(**)$ there are no nonzero analytic difference operators $L(x, \partial_x)$ and $\La(n, T)$
that satisfy one of the above conditions with $f(n) \equiv 0$ or $\theta(x) \equiv 0$.

 The assumption $(**)$ implies that the map $b \colon \B_1 \to \B_2$, given by $b(P(x, \partial_x)):= S(n, T)$ is a well defined algebra anti-isomorphism.  The algebra
\[
A_1 := b^{-1}(K_2) 
\]
consists of the bispectral  operators corresponding to $\psi(x,z)$ (i.e., 
differential operators in $x$ having the properties \eqref{bisp3}) and the algebra  

\[
\quad A_2 := b(K_1)
\]
consists of the bispectral  operators corresponding to $\psi(x,n)$, i.e.  
difference operators in $n$   having the properties \eqref{bisp3}).

Below we present the  discrete-discrete  version of the general bispectral problem which is   suitable in the set-up of  vector orthogonal polynomial sequences, which are eigenfunctions of difference operators.  

We are interested in the case when,  for any fixed $n$,  the   function $\psi(x,n)$  defining the map $b$ is a  polynomial in $x$.  We additionally assume that, for all $n$,  polynomials $\psi(x,n)$   are the eigenfunctions of a fixed differential operator in the variable $x$ and that, for any fixed $x$ the function $\psi(x,n)$  is an eigenfunction of a difference operator in $n$. We know that such a situation occurs in case of the classical discrete orthogonal polynomials.

Let $\B_1$ be the    algebra   spanned over $\Cset$ by  $x$,     $D$ and  $D^{-1}$. Needless to say that the commutation relations in $\B_i$  play a crucial role.  In the case of$\B_1$ they are

\[
[D, x] = D, \;\; [D^{-1}, x] = - D^{-1},\; \;  [D, D^{-1}] = 0.
\]

 In the  same way we define another algebra $\B_2$, using the operators   $T$,   its inverse $T^{-1}$   and the operator $n$ of multiplication by  the variable $n$.    In our case such a module is a  linear space of bivariate functions  $f(x,n),$ where $x$ is   a continuous   and $n$ is  a  discrete variables. For later use we also introduce $\Delta = D - 1$ and $\nabla = D^{-1} - 1$.   Let us put $L = x\nabla$. The above commutation relations yield

\[
[D, L]= - \Delta = 1 - D, \;\; [L, x] =L - x. 
\]

 In what follows we use the notation 
 
 \[
 (x)_k = x(x-1) \ldots (x-k+1)\;\; \text{for}\;\; k \in \Nset,   \;\; \text{and} \;\; (x)_0 =1.
 \]
 We notice that the  notation     $ (x)_k $ is quite often used with a different meaning but here we will use it only in the above sense.
 Let  $\psi(x,n):= S_n(x):=   (x)_n$. Set $L := x\Delta$   and $\La = T +n$.

 Obviously 

\[
L S(x,n)=  nS(x,n),\; \Lambda S(x,n) =  xS(x,n).
\]

In this way we can define the anti-automorphism  $b$ by

\begin{equation}    \label{iso}
\begin{cases}
b(x) =  T  + n\\
b(L)=   n\\  
b(\Delta) =   n T^{-1}. \\
\end{cases}
\end{equation}

%%%%%%%%%%%%%%%%%%%%%%%%%%%%%%%%%%%%%%%%%%%%%%%%%%%%%%%%%%%%%%%%%%%%%%%%%%%%%
%%%%%%%%%%%%%%%%%%%

In this subsection we, following \cite{BHY} recall  how to construct new bispectral operators from already  known ones. The method is quite general and does not depend on the specific form of the operators. First, we  remind the reader that for an operator $L \in \B$ it is said that $\ad_L$ acts     locally nilpotently on $\B$  when for  any element $a\in \B$ there exists $k\in \Nset$,  such that

\[
\ad^k_L(a)= 0.
\]

\medskip
The main  tool for constructing new bispectral operators in this paper will be the following simple observation made in \cite{BHY}. We formulate it in a form suitable for the  discrete VOP.
\bpr{2.1}   Let $\B_1, B_2$ be unital algebras with the properties described above.
Let $L \in \B_1$  and     $\ad_L : \B_1 \rightarrow \B_1$ be a locally nilpotent operator and let     $b\colon \B_1 \to \B_2$ be a bispectral involution. Suppose that, for any fixed $n$,  $e^L\psi(x,n)$ is a polynomial in $x$ of degree $n$. Define a new map $b^{'} :\B_1 \rightarrow \B_2$ via the new polynomial function $\psi'(x,n) := e^{\ad L} \psi(x,n)$. 
Then $b^{'}: \B_1\to \B_2$ is  a bispectral anti-involution.  
\epr

%%%%%%%%%%%%%%%%%%%%%%%%%%%%%%%%%%%%%%%%%%%%%%%%%%%%%%%%%%%%%%%%%%%%%%%%%%%%%
%%%%%%%%%%%%%%%%%%%

\section{Charlier type vector orthogonal polynomials}

 Here the algebras $\B_i$   are the ones defined in the previous section, i.e.  Let $P(X)$  be a polynomial of degree $d \geq 1$ without a free term.   We define the automorphism $\sigma:  \B_1 \rightarrow \B_1$ by

\[
\sigma  =  e^{\ad_{P(\Delta)}}.
\]
Let us compute explicitly its action  on the generators  $x$ and $\sigma(D)$.

\ble{aut-1}

The automorphism $\sigma$ acts on the generators  as
\[
\begin{cases}
\sigma(x)   =  x + P'(\Delta)D \\
\sigma(D) = D.
 \end{cases}
\]
\ele

\proof
Starting with the relation $[D, x]= D$ we
 prove by induction that for each $m$  

\[
[D^m, x] = mD^m.
\]
Again by induction we find that 

\beq
[\Delta^m, x]   = m \Delta^{m-1}D. \label{delta-comm}
\eeq
Really for $m=1$ it is obvious:
\[ 
 [\Delta, x]   = D.
\]
If \eqref{delta-comm} is verified for  $m=j-1$ we have for $j=m$

\begin{equation*}
\begin{aligned}
\Big[\Delta^m, x \Big]    & =  \Delta^{m} x - x \Delta^m = \Delta \Delta^{m-1} x - x \Delta^m   \\
&   =\Delta \{ x \Delta^{m-1} + (m-1)\Delta^{m-2} D\} - x \Delta^m \\
&  = \Big[\Delta, x \Big ]\Delta^{m-1}  +(m-1)\Delta^{m-2} D)  =    m \Delta^{m-1}D.
\end{aligned}
\end{equation*}

Hence

\[
 e^{\ad_{P(\Delta)}} (x)  = x + [P(\Delta), x] + \sum_{j=2}^{\infty} \ad_{P(\Delta)}^j (x) = x + P'(\Delta)D,
\]
as the rest of the terms vanish. The other formula is  obvious. \qed

A direct consequences of the lemma is 

\bco{aut-2}
The images of $L$  and $\Delta$ under the automorphism $\sigma$ are as follows:
\begin{eqnarray}
\nonumber \sigma(\Delta)& =& \Delta\\  
\sigma(L)& =& x\nabla + [P(\Delta), x \nabla  ]  =  x\nabla  -  P'(\Delta)\Delta. \label{sigma-L}
\end{eqnarray}
\eco

\proof 
The first identity is obvious. Here is the simple proof of the second one. We have
\[
[P(\Delta), x\nabla] =  (P(\Delta) x - x P(\Delta))\nabla =  P'(\Delta)D\nabla 
\]
as $\Delta$ and $\nabla$ commute. Then using that $D\nabla = - \Delta$ we obtain

\[
[P(\Delta), x\nabla] 
 =  - P'(\Delta)\Delta. 
\]
Hence

\[
\sigma(L) = L + [P(\Delta), x\nabla] = x\nabla - P'(\Delta)\Delta.
\]
\qed

Let us define the anti-involution  $b_1 = b( \sigma^{-1})$. Also we define  the difference operator 
\beq
\tilde{L} = \sigma^{-1}(L)  = x\nabla + P'(\Delta)\Delta. \label{d-1}
\eeq
Here and further we use that

\[
\sigma^{-1} = \sum_{j=0}^{\infty} \frac{(-\ad_{P(\Delta)})^j}{j!}
\]

  From \leref{aut} and \coref{aut-2} it follows almost immediately
that
\ble{inv}
The anti-involution $b_1$  acts as
\ele
\[ 
\begin{cases}
b_1(x)  = T + n + P'(nT^{-1})(nT^{-1} +1)\\
b_1(\tilde{L}) = n \\
b_1(\Delta) = n T^{-1}.
\end{cases}
\]
\proof

We have

\[ 
 b_1(x)  =   b(\sigma^{-1}(x)) = b(x + P'(\Delta)D)
 = T + n +P'(nT^{-1})(nT^{-1}  + 1). 
\]
Next
\[
b_1(\tilde{L}) =  b(\sigma^{-1}\circ \sigma(L))= b(L)   =    n  
\]
Finally 

\[
b_1(\Delta) = b(\Delta) = n T^{-1}.
\]
\qed

Let us define  the "wave function"   

\beq
C_n^P(x)=  e^{P(\Delta))}\psi(x,n)    = \sum_{j=0}^{\infty} \frac{P(\Delta)^j (x)_n}{j!}.
\label{wf-1}
\eeq

Notice that the operator  $\Delta$ reduces the degrees of the polynomials by one unit. The same is true for $P(\Delta)$ (we recall that $P(X) $  has no free term).   This shows that the  sum \eqref{wf-1} is finite   and for this reason $C_n^P(x)$  is a polynomial.

Let us write explicitly $P(\Delta)$ as

\[
P(\Delta) = \sum_{j=1}^{d} \beta_j\Delta^j.
\]
We will list the basic  properties of the polynomials $C_n^P(x)$ in terms of the polynomial $P$ in  the next theorem

\bth{mult-Ch}
The polynomials $C_n^P(x)$ have the following properties:

(i) They satisfy $d+2$-term recursion relation
\[
xC_n^P(x) =   C_{n+1}^P(x)    + n C_n^P(x)  
      +     \sum_{j=1}^{d} (n)_j (j \beta_j + (j+1)\beta_{j+1} ) C^P_{n-j}.
\]

(ii) They are eigenfunctions of the difference operator  $\tilde{L}$ \eqref{d-1}

\[
   \tilde{L} C_n^P(x)  = n C_n^P(x).
\]

 (iii) The following lowering operator holds
 
 \[
 \Delta C_n^P(x) = nC_{n-1}^P(x)
 \]

\ethe

\proof

(i)  From \leref{aut-1} we have that 

\[
xC_n^P(x) = \Big\{T +n + P'(nT^{-1})(nT^{-1} +1)\Big\}C^P_n 
\]
Let us work out the expression $E = P'(nT^{-1})(nT^{-1} +1)C^P_n$. 
We have 
\begin{eqnarray*}
E &=& \sum_{j=1}^d j \beta_j(nT^{-1})^jC^P_n  + \sum_{j=0}^d  + (j+1)\beta_{j+1}(nT^{-1})^jC^P_n\\
&=&\sum_{j=1}^d (j \beta_j +  (j+1)\beta_{j+1})(nT^{-1})^jC^P_n.
\end{eqnarray*}
 In this way we obtain

\begin{equation*}
\begin{split}
E &= \sum_{j=1}^{d}   \beta_j(nT^{-1})^j C^P_n  + \sum_{j=0}^d (n) \beta_{j+1}(nT^{-1})^{n-j}C^P_n =  \\
&= \sum_{j=1}^{d} (j \beta_j + (j+1)\beta_{j+1} ) (nT^{-1})^{n-j}C^P_n\\
&= \sum_{j=1}^{d} (n)_j (j \beta_j + (j+1)\beta_{j+1} ) C^P_{n-j}.
\end{split}
\end{equation*}
 
 % % % % % % % % % % % % % % % % % % % % % % % % % % %
 
 (ii) From the definitions of $\tilde{L}$ and $C_n^P(x)$ we obtain
 
 \[
 \tilde{L}C_n^P(x) = e^P L e^{-P}e^{P}\psi(x,n) = e^P n \psi(x,n)=n  C_n^P(x)  .
 \]

 % % % % % % % % % % % % % % % % % % % % % % % % % % % % % % % % % % % % % % % % % % %

(iii) follows directly from \leref{inv}.

 \section{Meixner  type vector orthogonal polynomials}\label{mei}

 We  use the notation of the previous section $x$, $\Delta$, $L = x\nabla$ to define an algebra $\B_1$ of discrete operators. It will be spanned by the operators $x$, $L$, $G =L\Delta + \beta \Delta$, where $\beta$   is a constant. They satisfy the following commutation relations

 \beq
 \begin{cases}
 [L, x]= - x -L\\
  [G, x]= L - x\Delta + \beta D   \\
  [G, L] = - G \\
    \ad_G^2 (x) = - 2G.
 \end{cases}
 \eeq
To simplify the notation we define the operator $R: = [G,x]$, which will be present both in the computations and in the final formulas. 
We use the  wave function from the previous section, namely $\psi(x,n) = (x)_n$. It has
the properties:

\[
\begin{cases}
L\psi(x,n)= n\psi(x,n)\\
x\psi(x,n)= \psi(x,n+1) +n\psi(x,n)\\
G\psi(x,n) =(n^2 + \beta n)    \psi(x,n-1).
\end{cases}
\]

We sum up these properties in terms of the following anti-involution $b$.

\[
\begin{cases}
b(L)= n \\
b(x)= T +n\\
b(\Delta) = nT^{-1}\\
b(G) = (n^2 + \beta n)T^{n-1}\\
b(R) = n(-n + \beta)T^{n-1} + \beta.
\end{cases}
\]
This gives our initial bispectral problem.

 Let $P(X)$  be a polynomial of degree $d \geq 1$   without a free term.     We define the automorphism $\sigma:  \B_1 \rightarrow \B_1$ by
 
 \[
 \sigma  =  e^{\ad_{P(G)}}.
 \]
 In the next lemma we compute it  on the generators. 
 
 \ble{aut}
 
 The automorphism $\sigma$ acts on the generators  as
 \[
 \begin{cases}
 \sigma(x)   =  x + RP' - P''(G)G    - P'^2(G)G\\
 \sigma(L) = L - P'(G)G\\
 \sigma(G) = G.
 \end{cases}
 \]
 \ele
 
 \proof
We need to compute $\ad^j_P$   the first few values of $j$ until it becomes $0$. First we prove by induction that for each $m$  
 
 \[
 [G^m, x] = \sum_{j=0}^{m-1} G^jRG^{m-1 -j}.
 \]
 Really we have for $m=1$ this is one of the  commutation relations in $\B_1$. Next assuming that the statement is true for $k$ we have
 \begin{eqnarray*}
 [G^{k+1}, x] &=& G G^kx - xG^{k+1} = G \Big[x G^k  +\sum_{j=0}^{k-1} G^jRG^{k-1 -j}\Big] - xG^{k+1}\\
 & = &G\sum_{j=0}^{k-1} G^jRG^{k -j} + 
 [Gx -x G] G^k  = \sum_{j=0}^{k} G^jRG^{k -j}.
 \end{eqnarray*}
   Next we use  that
   \[
   [G^j, R] = -2j G^j.
   \]
   This easily follows from the fact that 
   \[
   [G, R] = -2G.
   \]
 Hence we find
 \[
 [G^m, x] =   \sum_{j=0}^{m-1} (R-2j)G^{m-1} =mRG^{m-1} - m(m-1)G^{m-1}. 
 \]
 This shows that    
 \[
 \ad_P(x)=  RP'(G) - P''(G)G.
 \]

  Now we easily compute  $\ad^2_P (x)$:
 
 \[
 \ad^2_P (x) = [P(G), RP''(G)G] =- 2P'^2(G)G.
 \]
From the expressions for $\ad^j_P, j=0, 1, 2$ we obtain the first identity. 

The second identity    is an easy consequence from   

\[
[G^m , L] = -m G
\]
which can be proved by induction starting with $[G, L] = -G$.
The last  identities is obvious.   \qed
 
 Let us define the anti-involution  $b_1 = b( \sigma^{-1})$. 
 Let us define the operator $ \tilde{L}$  that will be our bispectral operator. Having in mind Meixner polynomials let $c \in \Cset, \;\; c \neq 0,1$.
 
 \beq
 \tilde{L}  = (1 - c)\sigma(L) =    (1- c)( x\nabla  - P'(G)G)   \label{d-ce}
 \eeq
 The next lemma computes the action of $b_1$ on the needed elements.

 \ble{inv-2}
 The anti-involution $b_1$  acts as

 \[ 
 \begin{cases}
 b_1(x)  =     T +n -  P'(n(n + \beta)T^{-1})[n(-n + \beta)T^{-1} + \beta]\\
  + P''((n(n + \beta)T^{-1})n(n + \beta)T^{-1}    - P'^2(n(n + \beta)T^{-1})n(n + \beta)T^{-1}\\
 b_1(L) =  n + P'(-nT^{-1})nT^{-1}\\
 b_1 (G)=  n(n +\beta)T^{-1} \\
 b_1(R)= n(-n + \beta)T^{-1}   +\beta.
 \end{cases}
 \]
  \ele
 \proof
 The last two identities are direct consequences of the definitions of $R$ and $G$ together with the formulas for $b$. The more involved  first identities follow from the last two. Really, we have
 
 \[ 
 b_1(x)  =   b(\sigma^{-1}x)
 = b(x) - b(P')b(R) + b(P''(G)G    - P'^2(G)G).
 \]
 after which we put the expressions for $b(G)$ an $b(R)$.
 Finally  
 \[
 b_1(\tilde{L}) =  b((1 - c)\sigma^{-1}(\sigma L)   = (1 - c)b(L) =  (1 - c) n. 
 \]
 \qed

 We come to the definition of the VOP, i.e.  the "wave function"   
 
 \beq
 M_n^P(x)=  e^{P(G))}\psi(x,n)    = \sum_{j=0}^{\infty} \frac{P(G)^j (x)_n}{j!}.
 \label{wf2}
 \eeq
 We assume that $P(G) =  \alpha G^m + \ldots$ with $\alpha \neq 0$.
 
 Notice that the operator  $G$ reduces the degrees of the polynomials by one unit. This shows that the  sum \eqref{wf2} is finite (we recall that $P(X) $  has no free term)  and for this reason $M_n^P(x)$  is a polynomial.

 In order to simplify the notations we introduce the following function:
 
 \[
 Q(G) =   [P''(G)   - P'^2(G)]G.
 \]
 The basic  properties of the polynomials $M_n^P(x)$ are listed in  the following theorem

 \begin{theorem} \label{mult-Mei}
 Let $\beta \notin \Nset$.  	Then the polynomials $M_n^P(x)$ have the following properties:
 	
 	(i) They satisfy the recursion relation
 	\begin{eqnarray*}
 	xM_n^P(x) &=&   M_{n+1}^P(x) + n M_n^P(x) +     \\
 	& + & \{ -  P'(n(n + \beta)T^{-1})n(-n + \beta)T^{-1} + Q(n(n+\beta)T^{-1})\}  M_n^P(x) . 
 		 \end{eqnarray*}
 	
 	(ii) They are eigenfunctions of the difference operator  $\tilde{L}$ \eqref{d-ce}
 	
 	\[
 	\tilde{L} M_n^P(x)  = (1 - c)n M_n^P(x).
 	\]
 	
 	(iii) The    operator $G$ acts on them as   lowering operator

 	\[
 	G M_n^P(x) = n(n + \beta)M^P_{n-1}(x).
 	\]

 \end{theorem}

 \proof  is similar to the proof \thref{mult-Ch} and follows easily from  \leref{inv-2} . Therefore it is omitted.
 \qed

\bre{Kr}
Notice that  the above polynomial system is well defined for all values of  $\beta$ but it is not always VOP. For example, when $\beta = - N, \; N \in \Nset$   and $d=1$ we come to Kravchuk  system of  orthogonal polynomials, which contains only finite number of members.  Similar situation occurs  when $d > 1$  This  is discussed in the next section.
\ere

 \section{Examples}

\bex{ex1}

For the definitions and properties  of  orthogonal polynomials we follow mainly \cite{KLS}.
Let us consider the  Charlier  polynomial sequence $\{C_n(x, a)\}$, defined by

\[
C_n(x;a) = \pFq{2}{0}{-n, -x}{-  }{ - \frac{1}{a}}, \; a > 0.
\]
These polynomials are orthogonal    with respect to the scalar product: 

\[
(C_m, C_n ) =   \sum_{x= 0}^{\infty}  \frac{a^x}{x!}C_m(x, a) C_n(x,a)   = a^{-n} e^an!\delta_{m n}. 
\]

Charlier polynomials  satisfy the following difference equation.

\[
- ny(x) = ay(x+1) - (x+a)y(x)+xy(x  - 1), \; y(x) =  C_n(x;a).
\]
      The above equation means that Charlier polynomials are eigenfunctions of the following difference operator 

\beq
aD - (x+a)I + x D^{-1} = a\Delta + x\nabla.  \label{ChO}
\eeq
Let us define the the normalized polynomials $q_n$, i.e which have a coefficient $1$ at the highest degree of $x$:

\[
q_n(x) = (-a)^n C_n(x;a).
\] 

They satisfy the 3-terms recursion relation

\[
xq_n(x) = q_{n+1}(x)  +(n+a)q_n(x) +an  q_{n+1}(x),   
\]
We will show that Charlier polynomils are a particular case of our construction from the previous section.

Let    $P(\Delta)  = - a\Delta$. We find that 

\[
\tilde{L} = a\Delta +   x \nabla.
\]
Define the polynomials  

\[
C_n^P  =e^{-a\Delta}(x)_n = \sum_{j=0}^n \frac{(-a)^j\Delta^j(x)_n}{j!}    =    \sum  \frac{(-a)^{j}(-n)_j(x)_{n-j}}{j!}.
\]
Now using the fact 

\[
\frac{(-n)_j}{j!} = (-1)^n\frac{(-n)_{n-j}}{(n-j)!}
\]
and changing the summation index $j \rightarrow k = n-j$ we obtain

\[
C_n^P=   (-a)^{n}\sum_{k=0}^n \frac{(-n)_k(x)_k}{(-a)^k k!}.
\]
We see that these are   the normed Charlier polynomials denoted above by $q_n(x)$. This example, together with the construction of Appell polynomials in \cite{Ho}, motivates  the name "Charlier-Appell polynomials" for  general  $P$.

\eex

\bex{M}
 In the second example we take the  algebra from section \ref{mei}.    Let us take $P(G) = \alpha G$, G = $(x\nabla + \beta)\Delta$. Then we have 

\[
\tilde{L} = (1-c)[\alpha x(\Delta + \nabla) +  \alpha\beta\Delta +   x\nabla].
\]
\eex
If we  take   the constant  $\alpha$  to be

\[
\alpha  =  \frac{c}{1-c} 
\]
 and $\beta$ to be different from a negative integer   we obtain exactly  the Meixner operator
     
     \[
\tilde{L} = c (x+\beta)\Delta +   x\nabla.
\]
 as given in \cite{KLS}.

 In case $\beta = -N$, $N \in \Nset$ we obtain Kravchuk polynomials $K_0, \ldots, K_N$, which form finite set of orthogonal polynomials.
 
\bex{M2}
Let us again use the   settings from section \ref{mei}. We present here the simplest new example. Let us take $P(G) =  \alpha G^2/2$. Then the new polynomials  

\[
M^P_n(x) = \sum_{j=0}^{\infty} \frac{(\alpha G^2)^j (x)_n}{j!}
\]
are eigenfunctions of the operator $\tilde{L} $  

\[
\tilde{L}M^P_n(x) =    (1- c)( x\nabla  -\alpha (x\nabla\Delta + \beta\Delta)^2)  M^P_n(x) = -nM^P_n(x)  
 \]
Choosing 

\[
\alpha  =  \frac{c}{ 2\beta (1-c)} 
\]
The operator $\tilde{L}$ becomes

\[
\tilde{L} =   x\nabla  + \frac{c}{2\beta}(x\nabla\Delta)^2. 
\]

The recurrence relation reads

\begin{eqnarray*}
xM^P_n &=& M^P_{n+1} + nM^P_n +   2\alpha n(n+\beta)2M_{n-1}^P    \\
 &+ &   \alpha n(n-1)(n  +\beta)(n + 1 -\beta)2M_{n-2}^P +    4\alpha^2          (-1) ^3 (n)_3(n -1+\beta)_3 M_{n-3}^P, 
\end{eqnarray*}

\eex

\bex{M-K} { \bf  Kravchuk-like polynomials.}
In this example we investigate the case when $\beta = - N, \; N \in \Nset$.  We take $P(G)  = G^2 + a G$.  We have 
$P'(G) = 2G +a$, $P''(G) = 2$,      $Q   = (2 -   (2G + a)^2)G$.  The   recursion relation reads

\begin{eqnarray*}
xM_n =M_{n+1} +n M_n  - n(-n -N) M_{n-1}  +
  \ldots - 4   (n)_3(n-N)_3 M_{n-3}.
\end{eqnarray*}
We see that the polynomials satisfy $5$-term recursion relation. However  the  coefficient at  $M_{n-3}$ is zero for $n = N$, thus violating the  condition of  P. Maroni's theorem   \cite{Ma}. This shows that  the vector orthogonality is valid only for the polynomials $M_n, n= 0, \ldots, N$.
The situation with the general polynomial $P(G)$ is similar.

\eex
%%%%%%%%%%%%%%%%%%%%%% References %%%%%%%%%%%%%%%%%%%%%%%%%%%%%%%%%%%%%%%

%%%%%%%%%%%%%%%%%%%%%%%%%%%%%%%%%%%%%%%%%%%%%%%%%%%%%%%%%%%%%%%%%%%%%%%%%%%%%%%
%%%%%%%%%%%%%%%%%%%%%%%%%%%%%%%%%%%%%%%%%%%%%%%%%%%%%%%%%%%%%%%%%%%%%%%%%%%%%%

\begin{thebibliography}{55}
%%%%%%%%%%%%%%%%%%%%%%%%%%%%%%%%%%%%%%%%%%%%%%

\bibitem{AlS} W. A. AI-Salam, {\em Characterization Theorems for Orthogonal Polynomials}, in: Orthogonal Polynomials: Theory and Practice
edited by P. Nevai,  with the assistance of M. E. H. Ismail, Proceedings of the NATO Advanced Study Institute on Orthogonal Polynomials and Their Applications
Colombus, Ohio, U.S.A. May 22 - June 3, 1989. Kluwer Academic publishers, Dordrecht / Boston / London.



\bibitem{ApKu}  A. I. Aptekarev, A. Kuijlaars, {\em Hermite-Pad\'e approximations and multiple orthogonal polynomial ensembles}, Uspekhi Mat. Nauk,
2011, Volume 66, Issue 6(402), 123–190

\bibitem{ACVA} J.~Arves\'u, J.~Coussement, W.~Van Assche, {\em Some discrete multiple orthogonal polynomials}, Journal of Computational and Applied Mathematics 153 (2003) 19–45.


   \bibitem{BHY} B. Bakalov, E. Horozov, and M. Yakimov, {\em{General methods for constructing
bispectral operators}}, Phys. Lett. A 222(1-2) (1996) 59--66.


\bibitem{BCZ}  Ben Cheikh, Y. and  A.~Zaghouani, {\em Some discrete d-orthogonal polynomial sets}, J. Comput. Appl. Math. 156, 253-263,
2003.


\bibitem{BCO}     Y. Ben Cheikh, A. Ouni {\em Some generalized hypergeometric d-orthogonal polynomial sets}, J. Math. Anal. Appl. 343 (2008) 464–478.


\bibitem{BDK}  P.~Bleher, S.~Delvaux, A.~B.~J. Kuijlaars, {\em Random matrix model with external source and a constrained vector equilibrium problem}, Comm. Pure Appl. Math., 64 (2011),  116–160.

 \bibitem{Bo} S.~Bochner, {\em \"Uber Sturm-Liouvillesche Polynomsysteme},  Math. Z. 29 (1929) 730--736.



 
 
 
 \bibitem{DG}  J.~J.~Duistermaat, and F.~A.~Gr\"unbaum,   {\em{ Differential equation in the spectral parameter}}, Commun. Math. Phys. 103 (1998) 177--240.
 
 
 
 \bibitem{DdlI} A.~J.~Dur\'an, Manuel D. de la Iglesia {\em Constructing Bispectral Orthogonal Polynomials from the Classical Discrete Families of Charlier, Meixner and Kravchuk}, Constructive Approximation
 February 2015, Volume 41, Issue 1, pp 49-91
 
 
 
 \bibitem{GVZ}   V.~X.~Genest,  L.~Vinet,  A.~Zhedanov,  {\em d-Orthogonal polynomials and su(2)}, Journal of Mathematical Analysis and Applications
 Volume 390, Issue 2, 15 June 2012, Pages 472-487.
 

 
 \bibitem{GHH}  F.~A.~Gr\"unbaum, L.~Haine, E.~Horozov,  {\em{ Some functions that generalize the Krall-Laguerre   polynomials}}, 
 { Jour. Comp. Appl. Math.}, 106(2), (1999) 271--297. 
 
 
 \bibitem{GY}  F.~A.~Gr\"unbaum, M.~Yakimov,  
 {\em{Discrete bispectral Darboux transformations from Jacobi operators}}, Pacific J. Math.   204(2) (2002),  395--431.
 
 
 \bibitem{vAssc} Walter Van Assche,  Difference Equations for Multiple Charlier and Meixner Polynomials, in:  Proceedings of the Sixth International Conference on Difference Equations Augsburg, Germany 2001,  New Progress in Difference Equations,   Edited by Saber Elaydi , Bernd Aulbach , and Gerasimos Ladas    CRC Press 2004,  Pages 549-557.
 
 \bibitem{Hil}  E.H. Hildebrandt, Systems of polynomials connected with the Charlier expansion and  the Pearson differential equation, Ann. Math. Statistics, 2(1931), 379-439.

\bibitem{Ho} E.~Horozov, Automorphisms of algebras and Bochner's property for vector orthogonal polynomials, 	arXiv:1512.03898 [math.CA].
 
 
 \bibitem{KLS}  R. Koekoek  P. A. Lesky 
 R. F. Swarttouw,   Hypergeometric  Orthogonal Polynomials  and Their q-Analogues, Springer Heidelberg Dordrecht London New York, 2010.
 Krall
 \bibitem{Kr}   H.~L.~Krall, On orthogonal polynomials satisfying a certain fourth order differential equation, The Pennsylvania State College Studies, No.6, The Pennsylvania State College, State College, PA, 1940.

\bibitem{Lan}  O.~E.~Lancaster, Orthogonal polynomials defined by difference equations, American Journal of Mathematics, 63(1941), 185-207.


 \bibitem{Lee}  D.W. Lee,  Difference equations for discrete classical multiple orthogonal polynomials Journal of Approximation Theory,  Volume 150, Issue 2, February 2008, Pages 132–152.

\bibitem{Les}  P.~Lesky, \"Uber Polynomsysteme, die Sturm-Liouvilleschen Differenzengleichungen gen\"ugen.
Mathematische Zeitschrift 78, 1962, 439–445.

\bibitem{Ma} P.~Maroni,  L'orthogonalit\'e et les r\'ecurrences de polyn\^omes d'ordre sup\'erieur \`a  deux,  Annales de la facult\'e des sciences de Toulouse (1989)  Volume: 10, Issue: 1, page 105-139.

\bibitem{NSU}  Nikiforov, A. F., Suslov, S. K., and Uvarov, V. B., Classical Orthogonal Polynomials of a Discrete Variable Springer,
New York, 1991.



\bibitem{ST} B.~Shapiro, M~Tater, {\em Asymptotic zero distribution of polynomial solutions to
	degenerate exactly-solvable equations}, Segovia, 2014.

\bibitem{VIs} J.~Van Iseghem,  {\em Vector orthogonal relations. Vector QD-algorithm}. Journal of Computational and Applied Mathematics
Volume 19, Issue 1, Supplement 1, July 1987, Pages 141–150.


\bibitem{VZh}  L.~Vinet, and A.~Zhedanov, {\em Automorphisms of the Heisenberg–Weyl algebra	and d-orthogonal polynomials}, JOURNAL OF MATHEMATICAL PHYSICS · MARCH 2009.



%%%%%%%%%%%%%%%%%%%%%%%%%%%%%%%%%
\end{thebibliography}
\end{document}